\documentclass[11pt]{amsart}
\title{Intersection properties of relations}
\author{Paolo Lipparini}
\address{Dipartimento di Matematica, Viale della Ricerca Scientifoca,
II Universit\`a di Roma (Tor Vergata),
 ROME 
ITALY}
\thanks{The author has received support from MPI and GNSAGA. We thank G. Cz\'edli for stimulating correspondence and discussions} 

\email{lipparin@axp.mat.uniroma2.it}
\urladdr{http://www.mat.uniroma2.it/\textasciitilde lipparin}
\newtheorem{Theorem}{Theorem}[section]

\newtheorem{prop}[Theorem]{Proposition}
\newtheorem{thm}[Theorem]{Theorem}

\newtheorem{cor}[Theorem]{Corollary}

\theoremstyle{definition}

\newcommand{\alg}{\mathbf} 

\begin{document}

\begin{abstract} 
We derive some equalities for relations on the algebra
 A, 
under the assumption that every subalgebra of 
A  $\times$  A
is congruence modular.
\end{abstract} 

\maketitle

\section{Notations}\label{notat} 

$ \alpha , \beta \dots$ denote \textit{congruences} on
some algebra
 ${\alg A}$;
$ \Theta, \Gamma$ are used
for \textit{tolerances}
(reflexive, symmetric and admissible relations),
 while we reserve the letters 
$R, S$ to denote reflexive
(but otherwise arbitrary)
 binary relations.
The words 
\textit{admissible} and
\textit{compatible} will be used with 
the same meaning.

Already
\cite{JoAU} realized 
the importance of dealing with
 reflexive admissible relations
(cf. \cite[Theorem 2.16]{JoAU}). \cite{Ts} presents a very clear
discussion of the interplay between congruence identities and 
identities involving  reflexive admissible relations.
Notice that our notation mainly comes from \cite{Ts},
and differs from \cite{JoAU}.

We shall write
$a R b$  to mean that $(a,b)\in R$, and we will
 use chains of the
above notation: for example, $a \Theta b \alpha c R d$ means
$(a,b)\in \Theta$, $(b,c) \in \alpha $ and $(c,d) \in R$.        

Intersection is sometimes denoted by juxtaposition; in particular
$ \alpha \beta $ denotes the meet of the congruences $ \alpha $ and $ \beta $.

 $R^*$ denotes the transitive closure of the binary relation $R$; in particular,
$\Theta^*$ is the smallest congruence which contains the 
tolerance $\Theta$.
$\overline{R}$ denotes the smallest compatible relation containing 
$R$ (where $R$ is a binary relation on some algebra which 
should be clear from the context).
In particular, $\overline{\Theta\cup\Gamma}$ is the smallest tolerance which contains the tolerances
$\Theta $ and  $\Gamma$. 

$R + S$ denotes 
$\bigcup _{n \in N}  \underbrace{R \circ S \circ R \circ S \dots} 
_{n\ factors}$ Thus, $R + S$ is the transitive closure
of $R \cup S$, and even the transitive closure
of $ R \circ S$, since
$R$ and  $S$ are supposed to be reflexive. In particular, if $ \alpha , \beta $
are congruences, $ \alpha + \beta $ is the join of $ \alpha $ 
and $ \beta $ in the lattice of congruences, while,
for $ \Theta$, $\Gamma$ tolerances, 
$ \Theta + \Gamma$
is the smallest \textit{congruence} which 
contains both $ \Theta$ and $\Gamma$. 
Notice that 
$ \Theta + \Gamma$
is far larger than the join of $ \Theta$ and $\Gamma$
\textit{in the lattice of tolerances}.

$R^-$ denotes the \textit{converse} of $R$, that is,
$aR^-b$ if and only if $bRa$. In particular,
$R+R^-$ is the smallest equivalence relation containing
$R$.  

$Cg(R)$ is the smallest congruence containing $R$.
Notice that if $R$ and $S$ are compatible  
then $R \circ S$ is compatible; and, by an induction, we get that
$R+S$ is compatible, too.
Thus, if $R$ is compatible, then
$Cg(R)=R+R^-$. In general, for $R$
not necessarily compatible, 
$Cg(R)=Cg(\overline{R})=
\overline{R}+\overline{R}^-$.

\section
{Intersection properties}

\begin{thm}
\label{subrelpiu} 
Suppose that $\alg A$ is an algebra such that 
every subalgebra of $\alg A\times \alg A$ 
generated by $4$ elements
satisfies 
$ \beta(\gamma \circ \delta \circ \gamma )\subseteq
\beta\gamma+\delta $, for all congruences $\beta,\gamma, \delta $
with $\delta\leq \beta $.

Then $\alg A$ satisfies 
\[
 \alpha (R + S)
\subseteq 
\alpha (\overline{R \cup S^-})+\alpha (\overline{R^- \cup S})=
\alpha (\overline{R \cup S})+\alpha (\overline{R^- \cup S^-})=
\alpha (Cg(R)+Cg(S))
\] 
for all reflexive relations $R$ and $S$ and every congruence  
$ \alpha $.
\end{thm}

\begin{prop}
\label{subrel} 
Under the hypothesis of
Theorem
\ref{subrelpiu}, 
 $\alg A$ satisfies 
\[
 \alpha (R \circ S)
\subseteq 
\alpha (\overline{R \cup S^-}) +
\alpha (\overline{R^- \cup S})
\] 
for all reflexive relations $R$ and $S$ and every congruence  
$ \alpha $.
\end{prop}

\begin{proof}
Let $a,c \in  A$ and
$(a,c)\in\alpha\cap(R \circ S)$. Thus
$a\alpha c$, and there
is $b\in  A$ such that
$a R b S c$.

Consider the subalgebra $\alg B$ of $\alg A\times \alg A$
generated by the four elements
$(a,a), (a,b), (c,b), (c,c)$.

First, observe that if $(x,y)\in \alg B$
then $(x , y) \in \overline{R \cup S^-}$, since all the generators of
$\alg B$ are in $ R \cup S^-$, and, by definition, 
$\overline{R \cup S^-}$
is compatible.

We have that $\big((a,a),(a,b) \big)$ and $\big ((c,b),(c,c) \big)$ belong
to $(0 \times 1)_{|\alg B}$,
$\big((a,b),(c,b) \big)$ belongs to $(\alpha \times 0)_{|\alg B}$
and
$\big((a,a),(c,c) \big)$ belongs to $(\alpha \times \alpha) _{|\alg B}$.

The above relations imply that
$\big ((a,a),(c,c) \big)$ belongs to
$$(\alpha \times \alpha) _{|\alg B}  
\cap
\Big((0 \times 1)_{|\alg B} \circ (\alpha \times 0)_{|\alg B} 
\circ (0 \times 1)_{|\alg B}
\Big)$$

Since $(\alpha \times 0)_{|\alg B} \leq
(\alpha \times \alpha) _{|\alg B}$, 
by the hypothesis of the Theorem,
\begin{multline*}  
(\alpha \times \alpha) _{|\alg B}  
\cap
\Big((0 \times 1)_{|\alg B} \circ (\alpha \times 0)_{|\alg B} 
\circ (0 \times 1)_{|\alg B}
\Big)
\subseteq  \\
\Big((\alpha \times \alpha) _{|\alg B} \cap 
(0 \times 1)_{|\alg B}\Big) + (\alpha \times 0)_{|\alg B}= 
(0 \times \alpha)_{|\alg B} + (\alpha \times 0)_{|\alg B}
\end{multline*}

In conclusion, $((a,a),(c,c))$ belongs to
$$(0 \times \alpha)_{|\alg B} + (\alpha \times 0)_{|\alg B}
$$

This implies that there is some $n$, and there are pairs
 $(x_i , y_i) \in {\alg B} $ 
 \makebox{$(0\leq i \leq n)$}
such that 
$$
(a,a)=(x_0, y_0)  \quad\quad (x_n, y_n )=(c,c)
$$
$$
(x_i, y_i) \equiv  (x_{i+1}, y_{i+1}) \mod (0 \times \alpha)_{|\alg B}\ \ \ \ \ \ {\rm for}\ i
\ {\rm even}  
$$
$$
(x_i, y_i) \equiv  (x_{i+1}, y_{i+1}) \mod (\alpha \times 0)_{|\alg B}\ \ \ \ \ \ {\rm for}\ i
\ {\rm odd}  
$$
In other words,
$$
a=x_0=y_0 \quad \quad c=x_n=y_n
$$
$$
x_i=x_{i+1},\ \ \ \ \ \ y_i \alpha y_{i+1} \ \ \ \ \ \ {\rm for}\ i
\ {\rm even}  
$$
$$
x_i \alpha x_{i+1},\ \ \ \ \ \ 
y_i=y_{i+1} \ \ \ \ \ \ {\rm for}\ i
\ {\rm odd}  
$$

In particular, $a=x_0=x_1 \alpha x_2=x_3 \alpha x_4 \dots$,
and $a=y_0 \alpha y_1=y_2 \alpha y_3=y_4 \dots$,
hence $x_i \alpha y_j$ for all $i$'s and $j$'s, since
$ \alpha $ is a congruence, and both  
$x_i$ and $y_j$ are congruent to $a$ modulo $ \alpha $.

Moreover, since $(x_i, y_i) \in {\alg B}$,
then 
$(x_i , y_i) \in \overline{R \cup S^-}$
 for all $i$'s,
by the remark made after the definition of ${\alg B}$.

Hence, for all $i$'s,  
$(x_i , y_i) \in \alpha (\overline{R \cup S^-})$,
and 
$(y_i , x_i) \in \Big(\alpha (\overline{R \cup S^-})\Big)^-=
\alpha (\overline{R^- \cup S})$.

In conclusion,
the sequence
$$
a=x_0=x_1 \ \ \  y_1=y_2 \ \ \ 
 x_2=x_3\ \ \  y_3=x_4 \ \ \ \dots \ \ \ x_n=y_n=c
$$
witnesses that $(a,c) \in 
\alpha (\overline{R \cup S^-})+
\alpha (\overline{R^- \cup S})$.
\end{proof}

\begin{cor}
\label{wtip} 
\cite{Lpwtip}
Under the hypothesis of
Theorem
\ref{subrelpiu}, 
 $\alg A$ satisfies 
\[
(wTIP)\quad\quad
 \alpha  \Theta^*
=
(\alpha  \Theta)^*
\] 
for every tolerance $\Theta$ and every congruence  
$ \alpha $.
\end{cor}

\begin{proof}
One inclusion is trivial.
By taking $R=S=\Theta$ in Proposition
\ref{subrel}, 
we get 
$
\alpha(\Theta \circ \Theta) \subseteq
(\alpha\Theta)^*
$.
The conclusion follows by 
induction: see 
\cite[Lemma 3.3]{Lpwtip}
 for details; actually, the argument
comes from \cite{CH} and \cite{CH1}.
\end{proof}


\begin{cor}
\label{rr} 
Under the Hypothesis of
Theorem
\ref{subrelpiu}, 
 $\alg A$ satisfies 
\[
\alpha(R+R^-) \subseteq 
\alpha (\overline{R}+\overline{R}^-)=
\alpha \overline{R}+\alpha \overline{R}^-=
\alpha Cg(R)
\] 
for every reflexive relation $R$.
\end{cor}

\begin{proof}
The first inclusion, as well as the inclusion
$\alpha (\overline{R}+\overline{R}^-)\supseteq 
\alpha \overline{R}+\alpha \overline{R}^-$
are trivial.

Since
$\overline{R}\circ \overline{R}^-$
is a tolerance, we can apply 
Corollary \ref{wtip} with  
$\overline{R}\circ \overline{R}^-$ in place of $\Theta$, thus getting
$\alpha (\overline{R}+\overline{R}^-)=
\alpha (\overline{R}\circ \overline{R}^-)^*=
\big(\alpha (\overline{R}\circ\overline{R}^-)\big)^*\subseteq 
\alpha \overline{R}+\alpha \overline{R}^-$,
where the last inclusion follows from
Proposition \ref{subrel}, with $\overline{R}^-$ in place of $S$, since 
$\overline{R}^{- -}=\overline{R}$,
and since 
$(\alpha \overline{R}+\alpha \overline{R}^-)^*=
\alpha \overline{R}+\alpha \overline{R}^-$.

$Cg(R)=\overline{R}+\overline{R}^-$
holds in every algebra, as mentioned 
at the end of Section \ref{notat},
hence
$ \alpha Cg(R)=\alpha (\overline{R}+\overline{R}^-)$.
\end{proof}

\begin{proof}
[Proof of Theorem \ref{subrelpiu}]
Since
$(\overline{R\cup S^-})^- = \overline{R^-\cup S}$,
 we can apply 
Corollary \ref{rr} with  
$\overline{R\cup S^-}$
in place of $R$,  getting
$\alpha(R + S) \subseteq 
\alpha\big((\overline{R\cup S^-}) + (\overline{R^-\cup S})\big)=
\alpha\big((\overline{R\cup S^-}) + (\overline{R\cup S^-})^-\big)=
\alpha(\overline{R\cup S^-}) + \alpha (\overline{R\cup S^-})^-=
\alpha(\overline{R\cup S^-}) + \alpha (\overline{R^-\cup S})$.

Since 
$\overline{R \cup S^-} \subseteq \overline{R}\circ \overline{S}^- $,
$\alpha(\overline{R\cup S^-})\subseteq
\alpha(\overline{R}\circ \overline{S}^-)\subseteq
\alpha(\overline{R\cup S})+\alpha(\overline{R^-\cup S^-})$,
by Proposition \ref{subrel} with 
$\overline{R} $ in place of $R$ and
$\overline{S}^-$ in place of $S$,
and since $\overline{\overline{R}\cup\overline{S}}=\overline{R \cup S}$.
Similarly,
$\alpha(\overline{R^-\cup S})
\subseteq\alpha(\overline{S}\circ \overline{R}^-)\subseteq
\alpha(\overline{R\cup S})+\alpha(\overline{R^-\cup S^-})$,
hence
$\alpha (\overline{R \cup S^-})+\alpha (\overline{R^- \cup S}) \subseteq 
\alpha (\overline{R \cup S})+\alpha (\overline{R^- \cup S^-})$.
By replacing $S$ with $S^-$ in the inclusion just obtained, we get
the reverse inclusion.

For the last identity, 
$\alpha (Cg(R)+Cg(S))=
\alpha (Cg(R \cup S))=
\alpha(\overline{R\cup S})+\alpha(\overline{R\cup S})^-=
\alpha(\overline{R\cup S})+\alpha(\overline{R^-\cup S^-})$
by the last identity in Corollary \ref{rr},
with $R \cup S$ in place of $R$. 
\end{proof}

\def\cprime{$'$} \def\cprime{$'$}
\providecommand{\bysame}{\leavevmode\hbox to3em{\hrulefill}\thinspace}
\providecommand{\MR}{\relax\ifhmode\unskip\space\fi MR }
\providecommand{\MRhref}[2]{%
  \href{http://www.ams.org/mathscinet-getitem?mr=#1}{#2}
}
\providecommand{\href}[2]{#2}

\end{document}